\documentclass[12pt]{amsart}

\usepackage{amssymb}

\usepackage{enumitem}

\usepackage{graphicx}

\makeatletter
\@namedef{subjclassname@2020}{%
  \textup{2020} Mathematics Subject Classification}
\makeatother

\usepackage[T1]{fontenc}

\newtheorem{theorem}{Theorem}[section]

\newtheorem{lemma}[theorem]{Lemma}

\theoremstyle{definition}

\numberwithin{equation}{section}

\begin{document}


\baselineskip=17pt


\title{On closed embeddings in ${\mathbb P}^{\mathbb N} \cup {\mathbb Q}^{\mathbb N}$ }
\author{El\.{z}bieta Pol, Roman Pol and Miros{\l}awa Re\'{n}ska}

\address{Institute of Mathematics\\University of Warsaw\\Banacha 2\\02-097 Warszawa, Poland}
\email{e.pol@mimuw.edu.pl}

\address{Institute of Mathematics\\ University of Warsaw\\Banacha 2\\ 02-097 Warszawa, Poland}
\email{r.pol@mimuw.edu.pl}

\address{Faculty of Mathematics and Information Science\\Warsaw University of Technology\\ 
Koszykowa 75\\00-662 Warszawa, Poland}
\email{miroslawa.renska@pw.edu.pl}


\begin{abstract}
Let  ${\mathbb P}$ and ${\mathbb Q}$ be the irrationals and the rationals in $[0,1]$. 
We prove that if a separable metrizable $X$ is a union of two disjoint 
$0$-dimensional sets $E$, $F$, 
$E$ is absolutely $G_{\delta }$ and $F$ is absolutely $F_{\sigma \delta }$, then
there is a closed embedding $h:X \to {\mathbb P}^{\mathbb N} \cup {\mathbb Q}^{\mathbb N}$ with 
$E= h^{-1} ({\mathbb P}^{\mathbb N} )$, $F= h^{-1} ({\mathbb Q}^{\mathbb N} )$.

We prove also that, for 
$H = \bigl\{ x \in I^{\mathbb N} : \forall k \ \exists l \ x (2^k 3^l ) = 0 \bigr\}$, 
whenever $A$ is an $F_{\sigma \delta }$-set in a compact 
one-dimensional metrizable space $X$, there is an embedding 
$h:X\to (\mathbb P \cup \{0 \} )^{\mathbb N} \cup {\mathbb Q}^{\mathbb N}$
such that $A = h^{-1} (H)$.

\end{abstract}

\maketitle
\section{Introduction}

We denote by ${\mathbb P}$ and ${\mathbb Q}$ the irrationals and the rationals in $I=[0,1]$, respectively, and 
${\mathbb P}^{\mathbb N}$, ${\mathbb Q}^{\mathbb N}$ are the countable products, the subspaces of the Hilbert cube 
$I^{\mathbb N}$.

All our spaces are metrizable separable, and our terminology follows  \cite{En} and \cite{Ku} .

In particular, $X$ is absolutely 
$G_{\delta }$ (absolutely 
$F_{\sigma \delta }$) if for some, equivalently - for any, embedding 
of $X$ in a completely metrizable space $Y$, the image of 
$X$ in $Y$ is a $G_{\delta }$-set ($F_{\sigma \delta }$-set, 
respectively); absolutely $F_{\sigma }$-sets are exactly the 
$\sigma $-compact sets, cf. \cite{En}, 4.5.7. 

The space ${\mathbb P}^{\mathbb N}$ is absolutely 
$G_{\delta }$ and ${\mathbb Q}^{\mathbb N}$ is absolutely 
$F_{\sigma \delta }$, cf. Comment 4.1.

The following result was proved in \cite{PPR}, Theorem 1,1.

\medskip

\begin{theorem} Let $\Gamma = \bigl\{ x\in {\mathbb Q}^{\mathbb N}: x(i) \in \{0,\ 1\}$ for 
all but finitely many $i \bigr\}$. Then for any disjoint zero-dimensional sets 
 $G$, $S$ in a compact space $X$, where $G$ is absolutely 
$G_{\delta }$ and 
$S$ is  $\sigma$-compact,  
there is an embedding $h:X\to I^{\mathbb N}$ such that 
$G = h^{-1} ({\mathbb P}^{\mathbb N})$ and 
$S = h^{-1} ({\Gamma})$.
\end{theorem}

\medskip

In this note we shall derive from Theorem 1.1 the following result.

\begin{theorem} Let $X$ be a union of two disjoint
zero-dimensional sets $E$, $F$,  where $E$ is absolutely 
$G_{\delta }$ and $F$ is absolutely 
$F_{\sigma \delta }$. 
Then there is an embedding 
$h:X\to {\mathbb P}^{\mathbb N} \cup {\mathbb Q}^{\mathbb N}$ 
onto a closed subspace such that 
$E = h^{-1} ({\mathbb P}^{\mathbb N})$ and 
$F = h^{-1} ({\mathbb Q}^{\mathbb N})$.
\end{theorem}

\medskip

In particular, ${\mathbb P}^{\mathbb N} \cup {\mathbb Q}^{\mathbb N}$, 
which is a union of a zero-dimensional absolutely 
$G_{\delta }$-set and a zero-dimensional absolutely 
$F_{\sigma \delta }$-set, contains topologically closed copies of every such space. 

We shall prove also another result concerning embeddings of one-dimensional 
absolutely $F_{\sigma \delta }$-spaces, using some key elements of the reasoning in 
\cite{PPR}, combined with Theorem 1.2.

Before we state this result, let us explain some background. Let 

\bigskip 

(1) $H_{kl} = \bigl\{ x \in I^{\mathbb N} : \ x (2^k 3^l ) = 0 \bigr\}$, 
$H = \bigcap\limits_k  \bigcup\limits_l H_{kl}$.

\bigskip

Then $H$ is a $F_{\sigma \delta }$-set in $I^{\mathbb N}$ and, by 
some standard arguments, one checks that for any $F_{\sigma \delta }$-set 
$A$ in a compact space $X$, there is an embedding 
$h:X\to {I}^{\mathbb N}$ with $A = h^{-1} (H)$, 
cf. \cite{Si}.

If $X$ is a compact zero-dimensional space and $A$ is a $F_{\sigma \delta }$-set 
in X, there is an embedding 
$h:X\to \{ 0, 1, {1\over 2}, {1\over 3}, \ldots \}^{\mathbb N}$ with 
$A = h^{-1} (H)$, cf. \cite{Ke}, sec. 23.A. 

We shall obtain a counterpart of these results for one-dimensional spaces.

\medskip

\begin{theorem}
For each  $F_{\sigma \delta }$-set 
$A$ in a compact one-dimensional space $X$, there is an embedding 
$h:X\to (\mathbb P \cup \{0 \} )^{\mathbb N} \cup {\mathbb Q}^{\mathbb N}$
such that $A = h^{-1} (H)$, 
cf. (1).
\end{theorem}

\medskip

In particular, 
$\bigl( (\mathbb P \cup \{0 \} )^{\mathbb N} \cup {\mathbb Q}^{\mathbb N} \bigr) \cap H$ 
is a one-dimensional absolutely $F_{\sigma \delta }$-set containing closed 
topological copies of any absolutely $F_{\sigma \delta }$-set of 
dimension $\leq 1 $, cf. Comment 4.1.

\section{Proof of Theorem 1.2.} 

\begin{proof}
Let $X = E \cup F$, where $E$, $F$ are disjoint
zero-dimensional sets, $E$ is absolutely 
$G_{\delta }$ and $F$ is absolutely 
$F_{\sigma \delta }$. 
Using \cite{Ku}, vol. II, \textsection 45, VII, Theorem 5, 
one can extend $X$ to a compact space $X^*$ with 
${\rm dim} X^* = {\rm dim} X \leq 1$
such that $F$ can be extended to a ${\sigma }$-compact 
zero-dimensional set $\widetilde F$ in $X^*$, disjoint from $E$. 

Indeed, $F$ being a 
$F_{\sigma }$-set in $X$, $F = \bigcup\limits_n F_n$, where 
$F_n$ is closed in $X$, and the cited theorem in \cite{Ku} provides a 
compact extension $X^*$ of $X$ with 
${\rm dim} X^* = {\rm dim} X $ such that for the closure 
$F_n^*$ of each $F_n$ in $X^*$, ${\rm dim} F_n^* = {\rm dim} F_n = 0$.

By absoluteness, $E$ is a $G_{\delta }$-set in $X^*$, 
and hence $\widetilde F = \bigcup\limits_n  F_n^* \setminus E$ 
has required properties.

Furthermore, we can enlarge the set $\widetilde F $ to a 
$\sigma $-compact zero-dimensional set $C$,

\smallskip\begin{enumerate}
\item[(1)] $\widetilde F \subset C \subset X^* \setminus E$,
\end{enumerate}
such that there are $G_{\delta }$-sets 
$G_1$, $G_2, \ldots $ in $X^*$ with

\smallskip\begin{enumerate}
\item[(2)] $X^* \setminus C =\bigcup\limits_i  G_i$, 
${\rm dim} (E \cup G_i ) = 0$.
\end{enumerate}

To get $C$ and $G_i$, we start from picking a base 
$B_1, B_2, \ldots $ of $X^*$ with boundaries 
${\rm Fr} B_i$ disjoint from the zero-dimensional set $E$. 
For $G_1 = X^* \setminus ( \widetilde F \cup \bigcup\limits_i {\rm Fr} B_i)$, 
we have $E \subset G_1$ and ${\rm dim} G_1 = 0$. 
Since ${\rm dim}\ {\rm Fr} B_i \leq 1$, one can choose for each $i$ a 
zero-dimensional $\sigma $-compact set $C_i \subset {\rm Fr} B_i $ with 
${\rm dim}({\rm Fr} B_i \setminus C_i ) \leq 0$, and let 
$C = \widetilde F \cup \bigcup\limits_i C_i$, 
$G_{i+1} = {\rm Fr} B_i \setminus C $, $i=1,2, \ldots $

We have then (1), the first part of (2), and since $E \subset G_1$ 
and $G_i$ is relatively closed in $E \cup G_i$ for $i\geq 2$, we get 
also the second part of (2).

Having defined $C$, we shall recall that, by absoluteness, $F$ is a 
$F_{\sigma \delta }$-set in $C$, and therefore, there are $\sigma $-compact 
sets $L_k$ such that

\smallskip\begin{enumerate}
\item[(3)] $L_1 = C \supset L_2 \supset L_3 \supset \ldots$, 
$F = \bigcap\limits_k L_k$, 
\end{enumerate}

and let

\smallskip\begin{enumerate}
\item[(4)] $L_k =  \bigcup\limits_i L_{ki}$, $L_{ki}$ compact.
\end{enumerate}  
Let ${\phi } : {\mathbb N} \times {\mathbb N } \to {\mathbb N }$ 
be a bijection and, cf. (4), (3),

\smallskip\begin{enumerate}
\item[(5)] $D_{\phi (k,i)} = L_{ki} \setminus L_{k+1}$
\end{enumerate}

The set  $L_{ki}$ being compact and disjoint from $E$, cf. (1), 
$D_{\phi (k,i)}$ is relatively closed in $E \cup D_{\phi (k,i)}$, and hence 
\smallskip\begin{enumerate}
\item[(6)] ${\rm dim}(E \cup D_i ) = 0$, for $i=1,2,\ldots $ 
\end{enumerate}

We are now in a position to apply Theorem 1.1 to get 
embeddings $u_i: X^* \to I^{\mathbb N }$ 
and $v_i: X^* \to I^{\mathbb N }$, $i=1,2,\ldots $, such 
that with $G_0 = \emptyset$, we have, cf. (2), (5) and (6),

\smallskip\begin{enumerate}
\item[(7)] $E \cup G_i = u_i^{-1} ({\mathbb P}^{\mathbb N})$,  
$C\subset u_i^{-1} ({\mathbb Q}^{\mathbb N})$, 
\end{enumerate}

\smallskip\begin{enumerate}
\item[(8)] $E \cup D_i = v_i^{-1} ({\mathbb P}^{\mathbb N})$,  
$C\setminus D_i \subset v_i^{-1} ({\mathbb Q}^{\mathbb N})$. 
\end{enumerate}

We shall show that the embedding 

\smallskip\begin{enumerate}
\item[(9)] $h=(u_0, v_1, u_1, v_2, \ldots ): X^* \to (I^{\mathbb N})^{\mathbb N}$
\end{enumerate}

 satisfies the condition

\smallskip\begin{enumerate}
\item[(10)] $E = h^{-1} \bigl(({\mathbb P}^{\mathbb N})^{\mathbb N}\bigr)$,  
$F = h^{-1} \bigl(({\mathbb Q}^{\mathbb N})^{\mathbb N}\bigr)$.
\end{enumerate}

Since $E \cup G_0 = E$, the first part of (10) follows from (7) and (8).

From (3), (5), (7) and (8) we have also 
$F \subset h^{-1} \bigl(({\mathbb Q}^{\mathbb N})^{\mathbb N}\bigr)$.
Let $x \not \in F$. If $x \not \in C$, then by (2), $x\in G_i$ for some $i \geq 1$, 
and by (7), $h(x) \not \in ({\mathbb Q}^{\mathbb N})^{\mathbb N}$.
Let $x \in C \setminus F$, and using (3) and (4) pick $k$, $i$ such that 
$x \in  L_{ki} \setminus L_{k+1}$. Then by (5), for $j= {\phi }(k,i)$, 
$x \in D_j$, and by (8) and (9), 
$h(x) \not \in ({\mathbb Q}^{\mathbb N})^{\mathbb N}$.

This ends the verification of (10).

To complete the proof of the theorem, it is enough to fix a bijection 
$\sigma : {\mathbb N} \to {\mathbb N} \times {\mathbb N}$ and 
to notice that for the homeomorphism 
$g:(I^{\mathbb N})^{\mathbb N} = I^{{\mathbb N} \times {\mathbb N}} \to  I^{\mathbb N}$, 
$g(x) (n) = x \bigl( \sigma (n) \bigr)$, we have 
$g\bigl(({\mathbb P}^{\mathbb N})^{\mathbb N}\bigr) = {\mathbb P}^{\mathbb N}$, and 
$g\bigl(({\mathbb Q}^{\mathbb N})^{\mathbb N}\bigr) = {\mathbb Q}^{\mathbb N}$.

\end{proof}

\section{Proof of theorem 1.3.}
A basic ingredient of our reasonings is the following lemma, which we shall derive from Lemma 2.1 in \cite{PPR} and Theorem 1.1.

\begin{lemma}
For each nonempty closed set $F$ in 
${\mathbb P}^{\mathbb N} \cup {\mathbb Q}^{\mathbb N}$ there is 
a continuous function 
$f:{\mathbb P}^{\mathbb N} \cup {\mathbb Q}^{\mathbb N} \to I$ 
such that $f({\mathbb P}^{\mathbb N} \setminus F ) \subset {\mathbb P}$, 
$f({\mathbb Q}^{\mathbb N}) \subset {\mathbb Q}$ and 
$F=f^{-1} (0)$.
\end{lemma}

\begin{proof}

Since ${\mathbb Q}^{\mathbb N}$ is a zero-dimensional 
$F_{\sigma }$-set in 
${\mathbb P}^{\mathbb N} \cup {\mathbb Q}^{\mathbb N}$, similarly as at 
the beginning of the proof of Theorem 1.2, one can use 
\cite{Ku}, vol. II, \textsection 45, VII, Theorem 5, 
to get an embedding 
$u:{\mathbb P}^{\mathbb N} \cup {\mathbb Q}^{\mathbb N} \to I^{\mathbb N}$ 
such that 
$u({\mathbb Q}^{\mathbb N} )$ is contained in a zero-dimensional  
${\sigma }$-compact set $T$ in 
$I^{\mathbb N} \setminus u({\mathbb P}^{\mathbb N} )$. 

Let $K=\overline{u(F)}$ be the closure of $u(F)$ in $I^{\mathbb N}$ 
(notice that 
$K \cap u({\mathbb P}^{\mathbb N} \cup {\mathbb Q}^{\mathbb N}) = u(F)$).

Let us split $T \setminus K$ into pairwise disjoint compact sets 
$C_1$, $C_2, \ldots $ with ${\rm diam }C_i \to 0$ 
(cf. the beginning of Sec. 4 in \cite{PPR}), and let 
$q:I^{\mathbb N} \to I^{\mathbb N} / \sim$ be the quotient 
map determined by the equivalence relation $\sim$, where 
the equivalence classes are $C_i$, $i=1,2,\ldots$ and the singletons 
$\{ x \}$, $x \not \in \bigcup\limits_i C_i = T \setminus K$. 

Then 
$q \circ u : {\mathbb P}^{\mathbb N} \cup {\mathbb Q}^{\mathbb N} \to I^{\mathbb N} / \sim$ embeds 
${\mathbb P}^{\mathbb N} \setminus F$ onto a 
$G_{\delta }$-set $G$, 
$S=q \circ u ({\mathbb Q}^{\mathbb N} \setminus F)$ is a countable set 
disjoint from $G$, and $A = q(K)$ is a compact set in 
$I^{\mathbb N} / \sim$ missing $G$. Applying Lemma 2.1 from \cite{PPR}, we 
get a continuous map $w:I^{\mathbb N} / \sim \ \to I$ with 
$w(G) \subset {\mathbb P}$, $w(S) \subset {\mathbb Q}$ and 
$A = w^{-1} (0)$. One readily checks that $f=w\circ q \circ u$ has 
required properties.

\end{proof}

\medskip

Having proved Lemma 3.1, one can now justify Theorem 1.3 easily. 

Since by Theorem 1.1, each one-dimensional compact space embeds in 
${\mathbb P}^{\mathbb N} \cup {\mathbb Q}^{\mathbb N}$, we can assume that 

\smallskip\begin{enumerate}
\item[(1)] $X \subset {\mathbb P}^{\mathbb N} \cup {\mathbb Q}^{\mathbb N}$.
\end{enumerate} 

The set $A \subset X$ being $F_{\sigma \delta }$ in $X$, 

\smallskip\begin{enumerate}
\item[(2)] $A = \bigcap\limits_k \bigcup\limits_l F_{kl}$, 
$F_{kl}$ closed in $X$.
\end{enumerate} 

We shall apply to each 
$F_{kl}\subset {\mathbb P}^{\mathbb N} \cup {\mathbb Q}^{\mathbb N}$, 
cf. (1), Lemma 3.1, getting continuous functions $f_{kl} : X \to I$ such that 

\smallskip\begin{enumerate}
\item[(3)] $f_{kl} \bigl( {\mathbb P}^{\mathbb N} \cap (X \setminus F_{kl})\bigr) \subset {\mathbb P}$, 
$f_{kl} ( {\mathbb Q}^{\mathbb N} \cap X ) \subset {\mathbb Q}$, 
$F_{kl} = f_{kl}^{-1} (0)$.
\end{enumerate} 

Let 

\smallskip\begin{enumerate}
\item[(4)] ${\mathbb M} = \{ n\in {\mathbb N}:\ n\not = 2^k 3^l,\ k,l = 1,2,\ldots \}$, 
$\sigma: {\mathbb M} \to {\mathbb N}$ bijection 
\end{enumerate} 

and let us define 
$h:X \to I^{\mathbb N}$ by, cf. (4), 

\smallskip\begin{enumerate}
\item[(5)] $h(x)(n) = \left\{
\begin{array}{cl}
 f_{kl}(x), &{\rm if} \ n= 2^k 3^l,\\
 x\bigl( \sigma (n) \bigr),  & {\rm if} \ n \in {\mathbb M}.
\end{array}
\right.$ 
\end{enumerate} 

The second part of (5) guarantees that $h$ is an embedding.

If $n= 2^k 3^l$, $h(x)(n) = f_{kl}(x)$, by (5), and hence, by (3), 
$h(x)(n) \in {\mathbb P}$ for 
$x\in {\mathbb P}^{\mathbb N} \setminus F_{kl}$, 
$f_{kl} (x) = 0$ for $x \in F_{kl}$, and 
$f_{kl} (x) \in {\mathbb Q} \setminus \{ 0\}$ 
for $x \in {\mathbb Q}^{\mathbb N} \setminus F_{kl}$.

It follows that 
$h(X) \subset ({\mathbb P} \cup \{ 0\} )^{\mathbb N} \cup {\mathbb Q}^{\mathbb N}$ and moreover, 
$F_{kl} = f_{kl}^{-1} (0) = h^{-1} (H_{kl})$, cf. (1) in Section 1. 
By  (2), we get $h^{-1} (H) = \bigcap\limits_k \bigcup\limits_l F_{kl} = A$.

\section{Comments.}

{\bf 4.1. Concerning ${\mathbb Q}^{\mathbb N}$.} The space 
${\mathbb Q}^{\mathbb N}$ (which is zero-dimensional and 
absolutely $F_{\sigma \delta }$) contains closed topological copies of all 
zero-dimensional  absolutely $F_{\sigma \delta }$-spaces. 

The following elegant reasoning to this effect was kindly communicated to us by Jan van Mill.

Let $X$ be a zero-dimensional absolutely $F_{\sigma \delta }$-space.
We can assume that $X \subset 2^{\mathbb N}$ and 
$X = \bigcap\limits_n F_n$, where $F_1 \supset F_2 \supset \ldots $ are 
$\sigma $-compact sets. Since $F_n \times 2^{\mathbb N} \times {\mathbb Q}$ 
is a zero-dimensional 
$\sigma $-compact nowhere countable and nowhere locally compact space, 
it is homeomorphic to $2^{\mathbb N} \times {\mathbb Q}$, by a classical 
Alexandroff-Urysohn theorem, cf. \cite{vM}, Corollary 5.2.
In effect, $F_n$ embeds onto a closed subspace of $2^{\mathbb N} \times {\mathbb Q}$.
It follows that $\prod\limits_{n} F_n$ embeds onto a closed subspace of
$(2^{\mathbb N} \times {\mathbb Q})^{\mathbb N}$ and so does $X$, 
which can be identified with the diagonal of $\prod\limits_{n} F_n$.
Finally, $2^{\mathbb N} \subset {\mathbb Q}^{\mathbb N}$, so 
$(2^{\mathbb N} \times {\mathbb Q})^{\mathbb N}$ embeds onto a closed subspace of
${\mathbb Q}^{\mathbb N}$.

The result follows also readily from a difficult theorem of van Engelen, 
providing a topological characterization of ${\mathbb Q}^{\mathbb N}$, 
cf. \cite{vE}, Section 3.

\medskip

{\bf 4.2. Concerning Theorem 1.3.} The result in Theorem 1.3 can be extended 
to higher Borel classes $F_{\sigma \delta \sigma}$, 
$F_{\sigma \delta \sigma \delta},\ldots$

To illustrate this point, let us consider the class of absolutely 
$F_{\sigma \delta \sigma}$-sets.

For each triple 
$(k,l,m) \in {\mathbb N} \times {\mathbb N} \times {\mathbb N}$ 
let us define, cf. (1) in Section 1
$H_{klm} = \bigl\{ x \in I^{\mathbb N} : \ x (2^k 3^l 5^m) = 0 \bigr\}$, 
and $H' =\bigcup\limits_m \bigcap\limits_k  \bigcup\limits_l H_{klm}$.

Then, with minor adjustments of the proof of Theorem 1.3, we get 
the assertion of this theorem for $F_{\sigma \delta \sigma}$-sets 
$A$ in compact one-dimensional spaces $X$, where $H$ is replaced by $H'$.

\end{document}